\documentclass[11pt]{article}
\usepackage{amsmath, epsfig, cite,colordvi,xcolor}
\usepackage{graphicx,subfigure}
\usepackage{amssymb}
\usepackage{amsfonts}
\usepackage{multirow}
\usepackage{array}
\usepackage{latexsym}
\usepackage{graphicx}
\usepackage{epstopdf}
\usepackage{tikz}
\usepackage{float}
\usepackage{cases}
\usepackage{mathrsfs}[12pt]
\usepackage{bm, amscd,  amsfonts, amssymb, cite,texdraw}
\usepackage{amsmath, epsfig, cite, color, colordvi}
\usepackage{amsmath,amsfonts,amssymb,amsthm}
\usepackage[colorlinks,urlcolor=purple,linkcolor=red,anchorcolor=green,citecolor=blue]{hyperref}
\usepackage{enumerate}
\usepackage{titlesec}
\usepackage{appendix}
\theoremstyle{plain}
\newtheorem{thm}{Theorem}[section]
\newtheorem{lem}[thm]{Lemma}

\numberwithin{equation}{section}

\def\hp{\mathcal{P}}

\def\hn{\mathcal{N}}
\def\qed{\hfill \rule{4pt}{7pt}}
\def\pf{\noindent {\it{Proof.} \hskip 2pt}}

\def\red{}

\begin{document}

\begin{center}

{\bf\Large The shifting method and generalized Tur\'{a}n number of matchings}\\[20pt]

{\large Jian Wang} \\[10pt]

Department of Mathematics\\
Taiyuan University of Technology, Taiyuan 030024, P.R. China\\[6pt]

Emails: wangjian01@tyut.edu.cn
\end{center}

\noindent {\bf Abstract.} Given two graphs $T$ and $F$, the maximum number of copies of $T$ in an $F$-free graph on $n$ vertices is called the generalized Tur\'{a}n number,  denoted by $ex(n,T,F)$. When $T=K_2$, it reduces to the classical Tur\'{a}n number $ex(n,F)$. Let $M_{k}$ be a matching with $k$ edges and $K^{*}_{s,t}$ a graph obtained from $K_{s,t}$ by replacing the part of size $s$ by a clique of the same size. In this paper, we show that for  any $s\geq 2$ \red{and $n\geq 2k+1$},
\[
ex(n,K_s,M_{k+1})=\max\left\{\binom{2k+1}{s}, \binom{k}{s}+(n-k)\binom{k}{s-1}\right\}.
\]
\red{For} any $s\geq 1$, $t\geq 2$ \red{and $n\geq 2k+1$},
\[
ex(n,K_{s,t}^*,M_{k+1})=\max\left\{\binom{2k+1}{s+t}\binom{s+t}{t}, \binom{k}{s}\binom{n-s}{t}+(n-k)\binom{k}{s+t-1}\binom{s+t-1}{t}\right\}.
\]
Moreover, we also study the bipartite case of the problem. Let $ex_{bip}(n,T,F)$ be the maximum possible number of copies of $T$ in an $F$-free bipartite graph with each part of size $n$. We prove that for any $s,t\geq 1$ \red{and $n\geq k$},
\[
ex_{bip}(n,K_{s,t},M_{k+1})=\left\{
\begin{aligned}
&\binom{k}{s}\binom{n}{t}+\binom{k}{t}\binom{n}{s}, & \quad s\neq t, \\
&\binom{k}{s}\binom{n}{s},&\quad s=t.
\end{aligned}
\right.
\]
Our proof is mainly based on the shifting method.


\allowdisplaybreaks

\section{Introduction}

Given a graph $T$ and  a family of graphs $\mathcal{F}$, the maximum number of copies of $T$ in an $\mathcal{F}$-free graph on $n$ vertices is called the {\it generalized Tur\'{a}n number}, denoted by $ex(n,T,\mathcal{F})$. When $T=K_2$, it reduces to the classical Tur\'{a}n number $ex(n,\mathcal{F})$. For a single graph $F$, we write $ex(n,T, F)$ instead of $ex(n,T,\{F\})$. \red{In \cite{zykov}, Zykov determined $ex(n,K_s,K_t)$ exactly.} \red{Let $P_k$ be the path on $k$ vertices and $\mathcal{C}_{\geq k}$ the family of all cycles with length at least $k$. In \cite{luo}, Luo gave the upper bounds on $ex(n,K_s,\mathcal{C}_{\geq k})$ and  $ex(n,K_s,P_k)$, which generalized the Erd\H{o}s-Gallai's Theorem on paths and cycles \cite{erdos2}.} Recently, the problem of estimating generalized Tur\'{a}n number has received a lot of attention, see \cite{alon,gyhori,gerbner1,gerbner2,gerbner3,gishboliner,jiema,boning}. We refer the reader to \cite{alon} for more background and motivation.

\red{Let $M_k$ be a matching with $k$ edges. In \cite{erdos2}, Erd\H{o}s and Gallai
proved the following theorem.
\begin{thm}\label{eg-thm}(Erd\H{o}s and Gallai \cite{erdos2})
For  any  $n\geq 2k+1$,
\[
ex(n,M_{k+1})=\max\left\{\binom{2k+1}{2}, \binom{k}{2}+k(n-k)\right\}.
\]
\end{thm}}

\red{Let $G_1$ and $G_2$ be two disjoint graphs. The {\it join} of two graphs, denoted by $G_1\vee G_2$, is defined as $V(G_1\cup G_2)=V(G_1)\cup V(G_2)$ and $E(G_1\cup G_2)=E(G_1)\cup E(G_2)\cup \{xy\colon x\in V(G_1), y\in V(G_2)\}$. We denote by $K_n$ and $E_n$ the complete graph on $n$ vertices and the empty graph on $n$ vertices, respectively. For the lower bound of Theorem \ref{eg-thm}, $K_{2k+1}$ and $K_k\vee E_{n-k}$ are $M_{k+1}$-free graphs with the required number of edges.}

In this paper, we determine the exact value of the generalized Tur\'{a}n number of matchings. Precisely, we prove the following two theorems.
\begin{thm}\label{main1}
\red{For  any $s\geq 2$ and $n\geq 2k+1$,}
\[
ex(n,K_s,M_{k+1})=\max\left\{\binom{2k+1}{s}, \binom{k}{s}+(n-k)\binom{k}{s-1}\right\}.
\]
\end{thm}

Let $K^{*}_{s,t}$ be a graph obtained from $K_{s,t}$ by replacing the part of size $s$ by a clique of the same size.

\begin{thm}\label{main1-2}
\red{For any $s\geq 1$, $t\geq 2$ and $n\geq 2k+1$,}
\[
ex(n,K_{s,t}^*,M_{k+1})=\max\left\{\binom{2k+1}{s+t}\binom{s+t}{t}, \binom{k}{s}\binom{n-s}{t}+(n-k)\binom{k}{s+t-1}\binom{s+t-1}{t}\right\}.
\]
\end{thm}

The proofs are mainly based on the shifting method, which has been used in \cite{akiyama} to give a short proof of the Erd\H{o}s-Gallai's Theorem on matchings.

We also study the bipartite case of the problem. Let $ex_{bip}(n,T,F)$ be the maximum possible number of copies of $T$ in a bipartite $F$-free graph with each part of equal size $n$. We prove the following theorem as well.
\begin{thm}\label{main2}
For any $s,t\geq 1$ and $n\geq k$, we have
\[
ex_{bip}(n,K_{s,t},M_{k+1})=\left\{
\begin{aligned}
&\binom{k}{s}\binom{n}{t}+\binom{k}{t}\binom{n}{s}, & \quad s\neq t, \\
&\binom{k}{s}\binom{n}{s},&\quad s=t.
\end{aligned}
\right.
\]
\end{thm}

{\noindent \bf Notations and outline.} Let $[n]$ denote the set $\{1,2,\ldots,n\}$. Let $G$ be a simple graph. By $V(G)$ and $E(G)$ we denote the vertex set and the edge set of $G$, respectively.  Let $e(G)$ be the number of edges of $G$ and $\hn(G,T)$ the number of copies of $T$ in $G$. For any $x\in V(G)$, we denote by $d_{G}(x)$ the number of neighbors of $x$ in $G$. If the graph $G$ is clear under the context, then we use $d(x)$ instead of $d_G(x)$.  By $\nu(G)$ we denote the number of edges in a maximum matching of $G$. For $S\subset V(G)$, the subgraph induced by $S$ is denoted by $G[S]$. For two disjoint subsets $S,T\subset V(G)$, the induced bipartite graph between $S$ and $T$ is denoted by $G[S,T]$.

The rest of the paper is organized as follows. In Section 2, we introduce the shifting operation on graphs and give some properties of this operation. In Section 3, we prove Theorems \ref{main1} and \ref{main1-2}. In Section 4, we prove Theorem \ref{main2}.

\section{The shifting operation}
Suppose a graph $G$ has vertex set $V(G)=[n]$ and edge set $E(G)=\{e_1,e_2,\ldots,e_m\}$. Here, edges in $E(G)$ are viewed as two-element subsets of $V(G)$. For $1\leq i<j\leq n$ and $e\in E(G)$, we define a {\it shifting} operation $S_{ij}$ on $e$ as follows:
\[
S_{ij}(e) =\left\{\begin{aligned}
&\left(e-\{j\}\right)\cup \{i\},\  if \ j\in e,\ i\notin e\mbox{ and }\left(e-\{j\}\right)\cup \{i\} \notin E(G),\\
&e,\ \mbox{otherwise}.
\end{aligned}\right.
\]
Define $S_{ij}(G)$ to be a graph on vertex set $V(G)$ with edge set $\{S_{ij}(e)\colon e\in E(G)\}$.

It is easy to see that $e(S_{ij}(G))=e(G)$.  In \cite{akiyama}, Akiyama and Frankl proved the following lemma.
\begin{lem}\label{matching}
Let $G$ be a graph on vertex set $[n]$. Then for any $1\leq i<j\leq n$,
\[
\nu(S_{ij}(G))\leq \nu(G).
\]
\end{lem}

We further prove that the shifting operation cannot reduce the number of the copies of $K_s$ and $K^*_{s,t}$.
\begin{lem}\label{clique}
Let $G$ be a graph on vertex set $[n]$. For  $1\leq i<j\leq n$,
\[
\hn(S_{ij}(G),K_s)\geq \hn(G,K_s) \mbox{\quad and \quad} \hn(S_{ij}(G),K^*_{s,t})\geq \hn(G,K^*_{s,t}).
\]
\end{lem}
\pf Firstly, we show that $\hn(S_{ij}(G),K_s)\geq \hn(G,K_s) $. Let $C$ be an $s$-element subset of $[n]$ that forms an $s$-clique in $G$. It is easy to check that $C$ also forms an $s$-clique in $S_{ij}(G)$ for the following three cases:

(i) $j\notin C$;

(ii) $j\in C$ and $i\in C$;

(iii) $j\in C$, $i\notin C$, but for all $x\in C\setminus \{j\}$, we have $\{x,i\}\in E(G)$.

\red{For the remaining case when} $j\in C$, $i\notin C$ and there exists some $x\in C$ such that $\{x,i\}\notin E(G)$. Then let
\[
A=\{x\in C\setminus \{j\}\colon \{x,i\} \notin E(G)\}
\]
and $B=C\setminus \{j\} \setminus A$. On one hand, since for any $x\in A$, $\{x,i\}\notin E(G)$,  $C'=\left(C\setminus \{j\}\right)\cup \{i\}$ does not form an $s$-clique in $G$. On the other hand, for any $x\in A$, we have $S_{ij}(\{x,j\})=\{x,i\}\in E(S_{ij}(G))$. And for any $x\in B$, we have $\{x,i\}\in E(S_{ij}(G))$. It follows that $C'$ forms an $s$-clique in $S_{ij}(G)$. Therefore, for any $s$-clique $C$ in $G$, there is an unique corresponding $s$-clique in $S_{ij}(G)$. Moreover, \red{there does not exist two $s$-cliques of $G$ that are shifted into one $s$-clique of $S_{ij}(G)$.} It follows that $\hn(S_{ij}(G),K_s)\geq \hn(G,K_s)$.

Secondly, we consider the number of copies of $K^*_{s,t}$. Let $C_1$ and $C_2$ be two disjoint subsets  of $[n]$ such that $G[C_1]$ forms an $s$-clique and $G[C_1,C_2] $ forms a copy of $K_{s,t}$. Then, each copy of $K^*_{s,t}$ in $G$ can be identified by such an ordered pair $(C_1,C_2)$. Let $C=C_1\cup C_2$. If $j\notin C$, then $(C_1,C_2)$ also forms a copy of $K^*_{s,t}$ in $S_{ij}(G)$. Therefore, in the rest of the proof, we assume that
 $j\in C$. Now the proof splits into the following two cases.

 {\bf Case 1.} $j\in C_1$. If $i\notin C$ and  for all $x\in C\setminus \{j\}$ we have $\{i,x\}\in E(G)$, then $(C_1,C_2)$ also forms a copy of $K^*_{s,t}$ in $S_{ij}(G)$. If $i\notin C$ and there exists some $x\in C\setminus \{j\}$ such that $\{i,x\}\notin E(G)$. Then $(C_1\setminus \{j\}\cup \{i\},C_2)$ is a copy of $K^*_{s,t}$ in $S_{ij}(G)$ but not a copy of $K^*_{s,t}$ in $G$. If $i\in C_1$, then for any $x\in C\setminus \{j\}$ we have $\{i,x\}\in E(G)$. It follows that $(C_1,C_2)$ is also a copy of $K^*_{s,t}$ in $S_{ij}(G)$. If $i\in C_2$ and for all $x\in C_2\setminus\{i\}$ we have $\{x,i\}\in E(G)$, then  $(C_1,C_2)$ also forms a copy of $K^*_{s,t}$ in $S_{ij}(G)$. If $i\in C_2$ and there exists some $x\in C_2\setminus \{i\}$ such that $\{i,x\}\notin E(G)$. Then, $(C_1\setminus\{j\}\cup \{i\}, C_2\setminus\{i\}\cup \{j\})$ is a copy of $K^*_{s,t}$ in $S_{ij}(G)$ but not a copy of $K^*_{s,t}$ in $G$.

 {\bf Case 2.} $j\in C_2$. If $i\notin C$  for all $x\in C_1$ we have $\{i,x\}\in E(G)$, then $(C_1,C_2)$ also forms a copy of $K^*_{s,t}$ in $S_{ij}(G)$. If $i\notin C$ and there exists some $x\in C_1$ such that $\{i,x\}\notin E(G)$. Then $(C_1,C_2\setminus \{j\}\cup \{i\})$ is a copy of $K^*_{s,t}$ in $S_{ij}(G)$ but not a copy of $K^*_{s,t}$ in $G$. If $i\in C$, then the neighborhood of $j$ is contained in the neighborhood of $i$ in this copy of $K^*_{s,t}$. It follows that $(C_1,C_2)$ also forms a copy of $K^*_{s,t}$ in $S_{ij}(G)$.

  Moreover, \red{it can be checked that there does not exist two  copies  of $K^*_{s,t}$ in $G$ that are shifted into one copy of  $K^*_{s,t}$ in $S_{ij}(G)$.} Combining all the cases, we conclude that $\hn(S_{ij}(G),K^*_{s,t})\geq \hn(G,K^*_{s,t})$. Thus, the lemma holds.
\qed
\section{The generalized Tur\'{a}n number of matchings}
\red{In this section, we determine the exact value of $ex(n,K_s,M_{k+1})$ and $ex(n,K^*_{s,t},M_{k+1})$ by characterizing all the shifted graphs with given matching number. The following lemma will be used in our proof, which is due to Bondy and Chv\'{a}tal \cite{bondy}.}
\begin{lem}\cite{bondy}\label{degree}
Let $G$ be a graph on $n$ vertices. If $\nu(G+uv)=k+1$ and $d(u)+d(v)\geq 2k+1$, then $\nu(G)=k+1$.
\end{lem}

\red{Let $G$ be a graph on vertex set $[n]$. We call $G$ {\it a shifted graph} if $S_{ij}(G)=G$ holds for all $i,j$ with $1\leq i<j\leq n$. If $G$ is a shifted graph, then for any $\{x,y\}\in E(G)$, $x'< x$ and $x'\neq y$, we always have $\{x',y\}\in E(G)$. Otherwise, we have $S_{x'x}(G)\neq G$, a contradiction. }

\red{For $k+1\leq \ell \leq 2k+1$, define graph $H(n,k,\ell)$ on vertex set $[n]$ as follows.  Let $A=[\ell]$, $B = [n]\setminus A$ and $C=[2k+1-\ell]\subset A$. The edge set of $H(n,k,\ell)$ consists of all edges between $B$ and $C$ together with
all edges in $A$. In the following lemma, we characterize all shifted graphs with given matching number.}

\begin{lem}\label{shifted-graph}
\red{Let $G$ be a shifted graph on vertex set $[n]$ with $\nu(G)=k$, then $G$ is a subgraph of $H(n,k,\ell)$ for some $k+1\leq \ell \leq 2k+1$.}
\end{lem}

\pf
Let $G'$ be an $M_{k+1}$-free graph on vertex set $[n]$ with maximum number of edges that containing $G$ as a subgraph. Then apply the shifting operation $S_{ij}$ to $G'$ for all $i,j$ with $1\leq i< j\leq n$. Finally, we obtain a graph $\tilde{G}$. Since $E(G)\subset E(G')$ and edges in $E(G)$ cannot be shifted, then $G$ is also a subgraph of $\tilde{G}$.

{\noindent\bf Claim 1.}  Vertex subset $[k+1]$ forms a clique in $\tilde{G}$.

\pf
Suppose to the contrary, there exist $x_1,x_2$ with $1\leq x_1<x_2\leq k+1$ such that $\{x_1,x_2\}\notin E(\tilde{G})$. Since $\nu(\tilde{G})= k$, it follows that there  exists an edge $\{y_1,y_2\}\in E(\tilde{G})$ such that $y_2>y_1\geq k$. Then by $x_1\leq k \leq y_1$ and $\{y_1,y_2\}\in E(\tilde{G})$, we have $\{x_1,y_2\}\in E(\tilde{G})$. By $x_2\leq k+1 \leq y_2$, we have $\{x_1,x_2\}\in E(\tilde{G})$, a contradiction. Thus, $\{1,2,\ldots,k+1\}$ forms a clique in $\tilde{G}$.\qed

Let $\ell$ be the maximum integer such that $[\ell]$ forms a clique in $\tilde{G}$. Let $U=[\ell]$ and  $U' = [n]\setminus U$. If $\ell\geq 2k+2$, then we shall obtain a matching of size $k+1$ in $\tilde{G}$, a contradiction. Thus, we have $k+1 \leq \ell\leq 2k+1$.

{\noindent\bf Claim 2.}  $U'$ forms an independent set \red{in} $\tilde{G}$.

\pf
Suppose to the contrary, there exist $x_1,x_2$ with $\ell+1\leq x_1<x_2\leq n$ such that $\{x_1,x_2\}\in E(\tilde{G})$. Then for any $x\in U$, since $x\leq \ell <x_1$ and $\ell+1<x_2$, then $\{x,\ell+1\}$ is an edge of $\tilde{G}$. It follows that $[\ell+1]$ forms a clique of $\tilde{G}$, a contradiction. Thus, the claim holds.\qed

{\noindent\bf Claim 3. }  For any vertex $y\in U'$, we have $d_{\tilde{G}}(y)\leq 2k-\ell+1$.

\pf
Suppose to the contrary that $d_{\tilde{G}}(y)\geq 2k-\ell+2$ for some $y\in U'$. Since $U$ is a maximum clique and $y\notin U$, it follows that there exists some $x\in U$, such that $\{x,y\} \notin E(\tilde{G})$. Since $\tilde{G}$ is the one with maximum number of edges, we have $\nu(\tilde{G}+xy)=k+1$. Since $d_{\tilde{G}}(x)\geq \ell-1$ and $d_{\tilde{G}}(y)\geq 2k-\ell+2$, then $d_{\tilde{G}}(x)+d_{\tilde{G}}(y)\geq 2k+1$. By Lemma \ref{degree}, it follows that $\nu(\tilde{G})=k+1$, a contradiction. Thus,  the claim holds.\qed

Combining all the claims, we conclude that $G$ is a subgraph of $H(n,k,\ell)$.
\qed

{\noindent \bf Proof of Theorem \ref{main1}.} Since $K_{2k+1}$ and $K_k\vee E_{n-k}$ are $M_{k+1}$-free graphs with the required number of $s$-cliques, we only need to prove the upper bound.
Let $G$ be an $M_{k+1}$-free graph on vertex set $[n]$ with the maximum number of $s$-cliques. \red{Since adding edges cannot reduce the number of $s$-cliques}, we assume that $G$ is the one with maximum number of edges subject to $\nu(G)\leq k$ and $\hn(G,K_s)$ is maximum. Clearly, we have $\nu(G)= k$. Otherwise, by adding one edge to $G$, we get a new graph $G'$ with more edges and $\nu(G')\leq k$, a contradiction. \red{By Lemmas \ref{matching} and  \ref{clique}, we can further assume $G$ is shifted. Then by Lemma \ref{shifted-graph}, we obtain that $G$ is a subgraph of $H(n,k,\ell)$ for some $k+1\leq \ell \leq 2k+1$.}

If $s> 2k+1$, since $G$ is $M_{k+1}$-free, it follows that $\hn(G,K_s)=0$. If $k+2\leq s\leq 2k+1$, we have
\[
\hn(G,K_s)\leq\hn(H(n,k,\ell),K_s) = \binom{\ell}{s} \leq \binom{2k+1}{s}.
\]

If $2\leq s\leq k+1$, then
\[
\hn(G,K_s)\leq\hn(H(n,k,\ell),K_s)= \binom{\ell}{s}+(n-\ell)\binom{2k-\ell+1}{s-1}.
\]
Let $f(\ell) =\binom{\ell}{s}+(n-\ell)\binom{2k-\ell+1}{s-1}$. By considering the second derivative, it is easy to check that $f(\ell)$ is a convex function. \red{Therefore since $k+1 \leq \ell\leq 2k+1$, we have}
\begin{align*}
\hn(G,K_s)&\leq \hn(H(n,k,\ell),K_s)\\
&\leq \max\left\{f(2k+1),f(k+1)\right\}\\
&= \max\left\{\binom{2k+1}{s},\binom{k}{s}+(n-k)\binom{k}{s-1}\right\}.
\end{align*}

Combining all the cases, we obtain that for $s\geq 2$ and $n\geq 2k+1$,
\[
ex(n,K_s,M_{k+1}) \leq \hn(G,K_s) \leq \max\left\{\binom{2k+1}{s},\binom{k}{s}+(n-k)\binom{k}{s-1}\right\}.
\]
 Thus, the theorem holds.
\qed

{\noindent \bf Proof of Theorem \ref{main1-2}.} \red{Since $K_{2k+1}$ and $K_k\vee E_{n-k}$ are $M_{k+1}$-free graphs with the required number of copies of $K^*_{s,t}$, we only need to prove the upper bound. Let $G$ be an $M_{k+1}$-free graph on vertex set $[n]$ with the maximum number of copies of $K^*_{s,t}$. We assume that $G$ is the one with maximum number of edges subject to $\nu(G)\leq k$ and $\hn(G,K^*_{s,t})$ is maximum.
By Lemmas \ref{matching} and  \ref{clique}, we can further assume $G$ is shifted. Then by Lemma \ref{shifted-graph}, we obtain that $G$ is a subgraph of $H(n,k,\ell)$ for some $k+1\leq \ell \leq 2k+1$.}

Let $\Omega_{\ell}(K^*_{s,t})$ be the set of all $K^*_{s,t}$'s in $H(n,k,\ell)$, i.e.,
\[
\Omega_{\ell}(K^*_{s,t}) = \left\{(C_1,C_2)\colon |C_1|=s,\ |C_2|=t \mbox{ and }(C_1,C_2) \mbox{ forms a copy of } K^*_{s,t} \mbox{ in }H(n,k,\ell)\right\}.
\]
Let $U=[\ell]$, $U_0=[2k+1-\ell]$ and $U'=[n]\setminus U$. Now we enumerate the copies of $K^*_{s,t}$ in $H(n,k,\ell)$ by classifying $\Omega_{\ell}(K^*_{s,t})$ into three classes as follows:
\[
\left\{
\begin{aligned}
\hp_1&=\{(C_1,C_2)\in \Omega_{\ell}(K^*_{s,t}) \colon C_1\subset U_0\};\\
\hp_2&=\{(C_1,C_2)\in \Omega_{\ell}(K^*_{s,t}) \colon C_1\cap U'\neq \emptyset \};\\
\hp_3&=\{(C_1,C_2)\in \Omega_{\ell}(K^*_{s,t}) \colon C_1\cap (U\setminus U_0)\neq \emptyset \}.
\end{aligned}
\right.
\]
For the first class, since there are $\binom{2k+1-\ell}{s}$ ways to choose $C_1$ and $\binom{n-s}{t}$ ways to choose $C_2$, it follows that
\[
|\hp_1|= f_1(\ell)=\binom{2k+1-\ell}{s}\binom{n-s}{t}.
\]
For the second class, since $U'$ is an independent set, there is exactly one vertex in $U'$ belonging to $C_1$ and all the other vertices in $C_1\cup C_2$ are contained in $U_0$. It follows that
\[
|\hp_2|= f_2(\ell)=(n-\ell)\binom{2k+1-\ell}{s-1}\binom{2k+1-\ell-s+1}{t}.
\]
For the third class, there are  $\binom{\ell}{s}-\binom{2k+1-\ell}{s}$ \red{choices} for $C_1$ and $\binom{\ell-s}{t}$ \red{choices} for $C_2$. Thus, we have
\[
|\hp_3|=f_3(\ell)=\left(\binom{\ell}{s}-\binom{2k+1-\ell}{s}\right)\binom{\ell-s}{t}.
\]
It is easy to check that $f_1(\ell),f_2(\ell)$ and $f_3(\ell)$ are all convex functions in $\ell$. Let
\[
f(\ell)=f_1(\ell)+f_2(\ell)+f_3(\ell).
\]
Then, $\hn(H(n,k,\ell),K^*_{s,t})=f(\ell)$ and $f(\ell)$  is a convex function in $\ell$. Thus, we have
\begin{align*}
\hn(G,K^*_{s,t})&\leq \hn(H(n,k,\ell),K^*_{s,t})\\[5pt]
&\leq \max\{f(2k+1),f(k+1)\}\\[5pt]
&=\max\left\{\binom{2k+1}{s}\binom{2k+1-s}{t},\binom{k}{s}\binom{n-s}{t}+(n-k)\binom{k}{s-1}\binom{k-s+1}{t}\right\}\\[5pt]
&=\max\left\{\binom{2k+1}{s+t}\binom{s+t}{t},\binom{k}{s}\binom{n-s}{t}+(n-k)\binom{k}{s+t-1}\binom{s+t-1}{t}\right\}.
\end{align*}
Thus, we complete the proof.
\qed

\section{The bipartite case}
\red{In this section, we determine the exact value of $ex_{bip}(n,K_{s,t},M_{k+1})$.
The following version of the K\"{o}nig-Hall
Theorem will be used in our proof.}

\begin{thm}(K\"{o}nig-Hall\cite{lovasz})
\red{Let $G$ be a bipartite graph with $\nu(G) = k$. Then there exists a subset $T$ of the
vertices with $|T| = k$, such that all edges of $G$ are incident to at least one vertex of $T$.}
\end{thm}
{\noindent \bf Proof of Theorem \ref{main2}.} \red{Let $X,Y$ be two vertex sets of size $n$. Let $G(X,Y)$ be an $M_{k+1}$-free bipartite graph with the maximum number of copies of $K_{s,t}$. We further assume that $G$ is the one with maximum number of edges subject to $\nu(G)\leq k$ and $\hn(G,K_{s,t})$ is maximum. Clearly, we have $\nu(G)= k$. Then, by K\"{o}nig-Hall Theorem, there exists a subset $T\subset X\cup Y$ with $|T|=k$, such that all edges of $G$ are incident to at least one vertex of $T$.}

\red{Let $X_1= X\cap T$, $Y_1 =Y\cap T$, $X_2= X\setminus T$ and $Y_1 =Y\setminus T$. Define $G^*$ to be a bipartite graph on vertex sets $X$ and $Y$ such that $G^*[X_1,Y]$ and $G^*[X,Y_1]$ are complete bipartite graphs and $G^*[X_2,Y_2]$ is an empty graph. Clearly, $G$ is a subgraph of $G^*$. It follows that $\hn(G,K_{s,t})\leq \hn(G^*,K_{s,t})$.}

\red{Let $(S,T)$ be an ordered pair such that $S\subset X$ and $T\subset Y$ with $|S|=s$ and $|T|=t$, and $G^*[S,T]$ is a complete bipartite graph. Clearly, each copy of $K_{s,t}$ in $G^*$ is identified by such an ordered pair. Since $G^*[X_2,Y_2]$ is an empty graph, it follows that at least one of $S\subset X_1$ and $T\subset Y_1$ holds. Let $|X_1|=x$. Since $|X_1|+|Y_1|=|T|=k$, then $|Y_1|=k-x$. Thus, we have
\[
\hn(G^*,K_{s,t}) = \binom{x}{s}\binom{n}{t}+\binom{n}{s}\binom{k-x}{t} -\binom{x}{s}\binom{k-x}{t}.
\]
Let
\[
f_{s,t}(x)= \binom{x}{s}\binom{n}{t}+\binom{n}{s}\binom{k-x}{t} -\binom{x}{s}\binom{k-x}{t}.
\]
It can be checked that $f_{s,t}(x)$ is a convex function. Thus, for $s=t$, we have
\[
\hn(G,K_{s,s})\leq \hn(G^*,K_{s,s})\leq \max\{f_{s,s}(0), f_{s,s}(k)\}= \binom{k}{s}\binom{n}{s}.
\]
Let $g(x) = f_{s,t}(x)+f_{t,s}(x)$. Since $g(x)$ is also a convex function, then for $s\neq t$, we have
\[
\hn(G,K_{s,t})\leq \hn(G^*,K_{s,t})\leq \max\left\{g(0), g(k)\right\}= \binom{k}{s}\binom{n}{t}+\binom{k}{t}\binom{n}{s}.
\]
Moreover, the complete bipartite graph $K_{k,n}$ implies the lower bound. Thus, the theorem holds.}
\qed

\noindent{\bf Acknowledgements.}
We thank the referees for their helpful comments. We also would like to thank Dr. B. Ning and Dr. X. Peng for sharing their knowledge on this topic. The research is supported by National Natural Science Foundation of China (No. 11701407) and Shanxi Province Science Foundation for Youths (No. 201801D221028).

\end{document}